\documentclass{amsart}
\usepackage{amsmath,amssymb}

\newcommand{\Z}{\ensuremath{\mathbf Z}}

\newtheorem{theorem}{Theorem}
\newcommand{\bt}{\begin{theorem}}
\newcommand{\et}{\end{theorem}}
\newtheorem{lemma}{Lemma}
\newcommand{\bl}{\begin{lemma}}
\newcommand{\el}{\end{lemma}}
\newcommand{\be}{\begin{eqnarray}}
\newcommand{\ee}{\end{eqnarray}}
\newcommand{\beq}{\begin{equation}}
\newcommand{\eeq}{\end{equation}}
\newcommand{\benum}{\begin{enumerate}}
\newcommand{\eenum}{\end{enumerate}}
\newcommand{\bsp}{\begin{split}}
\newcommand{\esp}{\end{split}}
\newcommand{\ba}{\begin{array}}
\newcommand{\ea}{\end{array}}

\begin{document}

\title{Maximal Sidon sets and matroids}

\author{J. A. Dias da Silva}
\address{Departamento de Matem\'{a}tica\\
Faculdade de Ci\^{e}ncias\\
Universidade de Lisboa\\
Campo Grande, Bloco C6\\
1749-016 Lisboa, Portugal}
\email{perdigao@hermite.cii.fc.ul.pt}
\author{Melvyn B. Nathanson}
\thanks{M.B.N. is supported
in part by grants from the NSA Mathematical Sciences Program and the
PSC-CUNY Research Award Program.}
\address{Department of Mathematics\\
Lehman College (CUNY)\\
Bronx, New York 10468, USA}
\email{melvyn.nathanson@lehman.cuny.edu}

\begin{abstract}
Let $X$ be a subset of an abelian group and
$a_1,\ldots,a_h,a'_1,\ldots,a'_h$ a sequence of $2h$ elements of
$X$ such that $a_1 + \cdots + a_h = a'_1 + \cdots + a'_h$.  The
set $X$ is a Sidon set of order $h$ if, after renumbering, $a_i =
a'_i$ for $i = 1,\ldots, h.$  For $k \leq h,$ the set $X$ is a
generalized Sidon set of order $(h,k)$, if, after renumbering,
$a_i = a'_i$ for $i = 1,\ldots, k.$  It is proved that if $X$ is a
generalized Sidon set of order $(2h-1,h-1),$ then the maximal
Sidon sets of order $h$ contained in $X$ have the same
cardinality.  Moreover, $X$ is a matroid where the independent
subsets of $X$ are the Sidon sets of order $h$.
\end{abstract}

\maketitle

\section{An extremal problem for Sidon sets}
Let $A$ be a subset of an abelian group $\Gamma$.   Two
$h$-tuples $(a_1,\ldots,a_h)$ and $(a'_1,\ldots,a'_h)$ of elements
of $A$ are called {\em equivalent}, denoted \[ (a_1,\ldots,a_h) \sim
(a'_1,\ldots,a'_h),
\]
if there is a permutation $\sigma$ of the set $\{1, \ldots,h\}$
such that $a'_i = a_{\sigma(i)}$ for $i=1,\ldots,h.$ If the
$h$-tuples $(a_1,\ldots,a_h)$ and $(a'_1, \ldots,a'_h)$ are
equivalent, then $a_1 +\cdots + a_h = a'_1 + \cdots + a'_h$.  We
write
\[ (a_1,\ldots,a_h) \not\sim (a'_1,\ldots,a'_h)
\]
if the $h$-tuples $(a_1,\ldots,a_h)$ and $(a'_1,\ldots,a'_h)$ are
not equivalent.

The {\em $h$-fold sumset of $A$}, denoted $hA$, is the set of all
elements of $\Gamma$ that can be written as the sum of $h$ elements
of $A$, with repetitions allowed. For every $x \in \Gamma,$ the {\em
representation function} $r_{A,h}(x)$ counts the number of
inequivalent representations of $x$ as a sum of $h$ elements of $A$,
that is, the number of equivalence classes of $h$-tuples
$(a_1,\ldots,a_h)$ such that $x = a_1+\cdots + a_h.$

The set $A$ is called a {\em Sidon set of order $h$} or a {\em
$B_h$-set} if every element of the sumset $hA$ has a unique
representation as the sum of $h$ elements of $A$, that is, if
$r_{A,h}(x) = 1$ for all $x \in hA$.  This means that if
$a_1,\ldots,a_h, a'_1,\ldots,a'_h \in A$ and
\[
a_1 + \cdots + a_h = a'_1 + \cdots + a'_h,
\]
then $(a_1,\ldots,a_h) \sim (a'_1,\ldots,a'_h).$ A {\em Sidon set}
is a Sidon set of order 2.

Let $X$ be a subset of the group $\Gamma,$ and denote by
$\mathcal{B}_h(X)$ the set of all finite $B_h$-sets contained in
$X$. Every set is a $B_1$-set, and
\[
\mathcal{B}_h(X) \subseteq \mathcal{B}_{h-1}(X) \subseteq \cdots
\subseteq \mathcal{B}_2(X) \subseteq \mathcal{B}_1(X).
\]
Moreover,
\[
\{a\} \in B_h(X) \qquad \text{for all $a \in X$ and $h \geq 1.$}
\]

If $A$ is a Sidon subset of order $h$ contained in the group
$\Gamma,$ then for every $x \in \Gamma$ the {\em translation} $x+A =
\{x+a:a \in A\}$ and the {\em reflection} $x-A = \{x-a:a \in A\}$
are also Sidon sets of order $h$.

A classical problem in combinatorial and additive number theory is
to determine the cardinality of the largest Sidon set of order $h$
contained in the interval of integers $\{1,2,\ldots,n\}.$ More
generally, if $X$ is a finite subset of the integers or of any
abelian group $\Gamma$, it is an open problem to estimate the
maximum cardinality of a Sidon set of order $h$ contained in $X$.
Every $B_h$-subset of a finite set $X$ is contained in a maximal
$B_h$-set, but there can be maximal Sidon sets of different
cardinalities contained in $X$. For example, in the interval
$\{1,2,3,4,5,6,7\}$, the sets $\{1, 3,6,7 \}$ and $\{1,2,5,7 \}$ are
the two maximal Sidon subsets of size 4, but there are also exactly
18 maximal Sidon subsets of size 3, for example, $\{1,3,4\}.$ Erd\H
os and Tur\' an~\cite{erdo-tura41} proved that the maximum size of a
Sidon set of order 2 contained in $\{1,2,\ldots,n\}$ is $n^{1/2} +
o\left(n^{1/2}\right).$   Ruzsa~\cite{ruzs98b} has constructed
maximal Sidon subsets of this interval of cardinality $\ll (n\log
n)^{1/3}.$ (See Martin and O'Bryant~\cite{mart-obry04} for
constructions of finite Sidon sets of integers and
O'Bryant~\cite{obry04} for a survey of the recent literature.)

The purpose of this paper is to describe a class of finite sets,
called $B_{2h-1,h-1}$-sets, in which all maximal Sidon sets of order
$h$ have the same cardinality.  Indeed, we shall prove that every
$B_{2h-1,h-1}$-set is a matroid in which the $B_h$-sets are the
independent sets. A maximal independent set in a matroid is called a
basis, and all bases in a matroid have the same cardinality.

\section{Generalized Sidon sets of order $(h,k)$}
Let $X$ be a subset of an abelian group $\Gamma.$   Let $h$ and $k$
be positive integers with $k \leq h.$ The set $X$ is called a {\em
generalized Sidon set of order $(h,k)$} or a {\em $B_{h,k}$-set} if,
whenever $a_1,\ldots,a_h, a'_1,\ldots,a'_h \in X$ and
\begin{equation} \label{ms:1}
a_1 + \cdots + a_h = a'_1 + \cdots + a'_h,
\end{equation}
there exist sets $I, I' \subseteq \{1, \ldots,h\}$ with $|I| = |I'|=
k$ and a one-to-one map $\tau: I' \rightarrow I$ such that
$a'_{i'} = a_{\tau(i')}$ for all $i' \in I'.$  Let $J = \{1,
\ldots,h\} \setminus I$ and $J' = \{1, \ldots,h\} \setminus I'$.
Then $|J| = |J'|= h - k$.  Since $X$ is a subset of the group
$\Gamma,$ it follows by subtraction that
\begin{equation} \label{ms:2}
\sum_{j\in J} a_j = \sum_{j'\in J'} a_{j'}.
\end{equation}
The Sidon sets of order $h$ are precisely the $B_{h,h}$-sets.

A simple example of a $B_{3,1}$-set that is not a $B_3$-set is
$\{1,2,3\}$.  Indeed,
\[
\{1,2,3\} \in \mathcal{B}_{2h-1,1}(\Z)\setminus
\mathcal{B}_{2h-1,2}(\Z)
\]
for every $h \geq 2.$ Another example of a $B_{3,1}$-set is $\{1,
14, 19, 20, 25,38\}$.

Let $\mathcal{B}_{h,k} = \mathcal{B}_{h,k}(\Gamma)$ denote the set
of all finite $B_{h,k}$-sets in $\Gamma$. It follows
from~(\ref{ms:1}) and~(\ref{ms:2}) that if $A$ is a $B_{(h,k)}$-set
and also a $B_{h-k}$-set, then $A$ is a $B_{h}$-set.  Conversely, if
$A$ is a $B_{h}$-set, then $A$ is both a $B_{h,k}$-set and a
$B_{h-k}$-set.  Therefore,
\begin{equation}  \label{ms:iden1}
\mathcal{B}_h = \mathcal{B}_{(h,k)} \cap \mathcal{B}_{h-k}
\end{equation}
for $k = 1,\ldots, h.$  In particular, for $h \geq 2$ we have
\begin{equation}   \label{ms:iden2}
\mathcal{B}_{2h-1} = \mathcal{B}_{(2h-1,h-1)} \cap \mathcal{B}_{h}.
\end{equation}
Thus, if $A$ is a $B_h$-subset of a $B_{2h-1,h-1}$-set,
then $A$ is a $B_{2h-1}$-set.

Similarly, if $k$ and $\ell$ are positive integers and $k+\ell \leq
h,$ then
\begin{equation}  \label{ms:iden3}
\mathcal{B}_{h,k} = \mathcal{B}_{h,\ell} \cap \mathcal{B}_{h-\ell,k-\ell}
\end{equation}
for $0 \leq \ell < k \leq h.$
It follows that
\begin{equation}  \label{ms:iden4}
\mathcal{B}_{2h-1,h-1} \subseteq \mathcal{B}_{2h-k,h-k}
\end{equation}
for $1 \leq k \leq h-1.$

Let $X$ be a subset of an abelian group $\Gamma$.  We have
\begin{equation}  \label{ms:iden5}
\mathcal{B}_{h,h}(X) \subseteq \cdots \subseteq
\mathcal{B}_{h,k+1}(X) \subseteq \mathcal{B}_{h,k}(X) \subseteq
\cdots \subseteq \mathcal{B}_{h,1}(X)
\end{equation}
for $k = 1,\ldots,h-1.$ In the group \Z\ of integers, if $g> h,$
then every finite subset of the set $\{g^i:i=1,2,3,\ldots\}$ is a
$B_h$-set, and so $B_{h,k}$-sets exist for all $h \geq 1$ and $k =
1,\ldots,h.$   However, not all of the set inclusions
in~(\ref{ms:iden5}) are proper.

\bt Let $X$ be a subset of an abelian group $\Gamma$. If $h \geq 2$ and $k \geq h/2,$ then $\mathcal{B}_h(X) =
\mathcal{B}_{h,k}(X).$ \et

\begin{proof}
It suffices to show that $\mathcal{B}_{h,k+1}(X) =
\mathcal{B}_{h,k}(X)$ if $h/2 \leq k \leq h-1.$

If $\mathcal{B}_{h,k+1}(X) \neq \mathcal{B}_{h,k}(X)$, then there is a set $A \in \mathcal{B}_{h,k}(X) \setminus \mathcal{B}_{h,k+1}(X)$ and there are elements $a_1,\ldots,a_{h-k},a'_1,\ldots,a'_{h-k}\in A$ such that
\beq  \label{ms:ah-k}   a_1+\cdots + a_{h-k} = a'_1 +\cdots +
a'_{h-k} \eeq and
\[
\{a_1,\ldots,a_{h-k}\} \cap \{a'_1,\ldots,a'_{h-k}\} = \emptyset.
\]
The inequality  $h/2 \leq k \leq h-1$ implies that $1 \leq h-k \leq k$.  By the
division algorithm,
\[
h = q(h-k)+r,
\]
where $q \geq 1$ and  $0 \leq r < h-k.$  It follows from
~(\ref{ms:ah-k}) that
\[
qa_1+\cdots + qa_{h-k} + ra^\ast = qa'_1 +\cdots + qa'_{h-k} +
ra^\ast
\]
for any $a^\ast \in A$. Each side of this equation is a sum of $h$
elements of $A$, but the two sides have only $r <h-k \leq k$ common
summands. This is impossible if $A \in \mathcal{B}_{h,k}(X)$, and so
$\mathcal{B}_{h,k+1}(X) = \mathcal{B}_{h,k}(X)$. This completes the
proof.
\end{proof}

Dias da Silva and Nathanson~\cite{dias-nath05} have constructed
nontrivial generalized Sidon sets of order $(2h-1,h-1)$ for all $h
\geq 2.$

We remark that if $\mathcal{B}_{h,k}(\Z) \setminus
\mathcal{B}_{h,k+1}(\Z) \neq \emptyset$, then
$\mathcal{B}_{h,k}(\Z) \setminus \mathcal{B}_{h,k+1}(\Z)$ contains
arbitrarily large finite sets of integers.

\bt Let $1 \leq k < h/2$.  If $A$ is a finite set of integers in
$\mathcal{B}_{h,k} \setminus \mathcal{B}_{h,k+1}$ and $b
> h\max(A),$ then $A \cup \{ b \} \in \mathcal{B}_{h,k}
\setminus \mathcal{B}_{h,k+1}.$
\et

\begin{proof}
Let $A^\ast = A \cup \{ b \}.$ Let $0 \leq r \leq s \leq h$ and let
$\{a_i\}_{i=1}^{h-r}$ and  $\{a'_i\}_{i=1}^{h-s}$ be subsets of $A$
such that
\[
rb+ \sum_{i=1}^{h-r} a_i = sb + \sum_{i=1}^{h-s} a'_i.
\]
We must show that at least $k$ summands on the left are the same as
$k$ summands on the right.  If $0 \leq r < s \leq h,$ then
\[
rb+ \sum_{i=1}^{h-r} a_i < (r+1)b \leq sb \leq  sb +
\sum_{i=1}^{h-s} a'_i,
\]
which is absurd.  Therefore, $r = s$ and
\[
\sum_{i=1}^{h-r} a_i = \sum_{i=1}^{h-r} a'_i.
\]
If $r \geq k,$ we are done.  If $r < k,$ then $A \in
\mathcal{B}_{h,k} \subseteq \mathcal{B}_{h-r,k-r}$ implies that
$k-r$ summands on the left are the same as $k-r$ summands on the
right, and so $A \cup \{ b \} \in \mathcal{B}_{h,k}.$

If $A \not\in \mathcal{B}_{h,k+1}(\Z),$ then $A \cup \{ b \} \not\in \mathcal{B}_{h,k+1}(\Z).$
This completes the proof.
\end{proof}

\section{Maximal Sidon sets of order $h$}

Let $X$ be a subset of an abelian group $\Gamma.$  A {\em double
representation of length $\ell$ in $X$} is a sequence $a_1, a_2,
\ldots, a_{\ell}, a'_1, a'_2, \ldots, a'_{\ell}$ of $2\ell$ not
necessarily distinct elements of $X$ such that \beq
\label{ms:doublerep} a_1 + a_2 + \cdots + a_{\ell} = a'_1 + a'_2 +
\cdots + a'_{\ell} \eeq and \[ (a_1,\ldots,a_{\ell}) \not\sim
(a'_1,\ldots,a'_{\ell}).
\]
There exists a double representation of length $h$ in the set $X$
if and only if $X$ is not a $B_h$-set.

The double representation~(\ref{ms:doublerep}) is called {\em
proper} if
\[
\{ a_1,a_2,\ldots,a_{\ell} \} \cap  \{ a'_1,a'_2,\ldots,a'_{\ell} \}
= \emptyset.
\]
If~(\ref{ms:doublerep}) is a double representation of length $\ell,$
then we can cancel elements that appear on both sides of the
equation, and obtain a unique proper double representation of length
$\ell',$ where $1 \leq \ell' \leq \ell.$

\bl  \label{ms:lemma:proper-h}
Let $h \geq 2$ and let $X$ be a finite
$B_{2h-1,h-1}$-subset of an abelian group $\Gamma.$
If
\[
a_1 + a_2 + \cdots + a_{\ell} = a'_1 + a'_2 + \cdots + a'_{\ell}
\]
is a proper double representation of length $\ell \leq 2h-1,$ then
$\ell = h$. \el

\begin{proof}
By~(\ref{ms:iden4}) we have
\[
\mathcal{B}_{2h-1,h-1} \subseteq \mathcal{B}_{2h-k,h-k}
\]
for $k = 1,2,\ldots, h-1.$  If $h+1 \leq \ell \leq 2h-1,$ then
\[
\ell = 2h-k,
\]
where
\[
1 \leq k \leq h-1.
\]
Since
\[
X \in \mathcal{B}_{2h-k,h-k}
\]
and $h-k \geq 1,$ it follows that $a'_i = a_j$ for some $i,j \in
\{1,\ldots, \ell\},$ which contradicts the hypothesis that the
double representation is proper.  Therefore, $\ell \leq h.$

Suppose that $\ell \leq h-1.$  By the division algorithm, there exist integers $q$ and $r$ such that
\[
2h-1 = q\ell + r
\]
and
\[
0 \leq r \leq \ell - 1 \leq h-2.
\]
Then
\[
qa_1 + qa_2 + \cdots + qa_{\ell} = qa'_1 + qa'_2 + \cdots +
qa'_{\ell}
\]
is a proper double representation of length $q\ell,$ where
\[
h+1 \leq q\ell = 2h-1-r \leq 2h-1,
\]
which is impossible.  Therefore, $\ell = h.$  This completes the
proof.
\end{proof}

\bl  \label{ms:lemma:proper-one}
Let $h \geq 2,$ let $X$ be a finite
$B_{2h-1,h-1}$-subset of an abelian group $\Gamma,$ and let $A$ be a
{\em maximal} $B_h$-subset of $X$.
For every
\[
x \in X \setminus A,
\]
there is exactly one proper double representation with elements
in  $A \cup \{ x \}$ and of length at most $2h-1.$
\el

\begin{proof}
Since $A$ is a maximal $B_h$-set contained in $X$, it follows that
$A\cup \{x\}$ is not a $B_h$-set, and so there exists a double
representation of the form
\[
ux + a_1+\cdots + a_{h-u} = vx + a'_1 + a'_2 + \cdots + a'_{h-v}.
\]
Let $u \geq v.$  Subtracting equal elements that appear on both sides of this
equation and renumbering the elements that remain in the equation,
we obtain a proper double representation of length $\ell \leq h$.
By Lemma~\ref{ms:lemma:proper-h}, we must have $\ell = h$, and so $v = 0$, there is no cancelation,
and the proper double representation is be of the form
\beq \label{ms:uniqueu}
ux + a_1+\cdots + a_{h-u} = a'_1 + a'_2 + \cdots + a'_h.
\eeq

Suppose that $w \geq 1$ and
\beq \label{ms:uniquew}
wx + b_1+\cdots + b_{h-w} = b'_1 + b'_2 + \cdots + b'_h
\eeq
is also a proper double
representation of length $h$ in $A \cup \{ x \}.$ Adding
equations~(\ref{ms:uniqueu}) and~(\ref{ms:uniquew}), we obtain
\[
ux + a_1+\cdots + a_{h-u}+ b'_1 + b'_2 + \cdots + b'_h = wx + a'_1 +
a'_2 + \cdots + a'_h + b_1+\cdots + b_{h-w}.
\]
If $u = w,$ we obtain the relation
\[
a_1+\cdots + a_{h-u}+ b'_1 + b'_2 + \cdots + b'_h = a'_1 + a'_2 +
\cdots + a'_h + b_1+\cdots + b_{h-u},
\]
where all of the summands belong to the $B_h$-set $A$.  It follows
that every term on the left appears on the right, and conversely.
Since $a_j \neq a'_i$ for all $i$ and $j$, we must have a
bijection between the sets $\{ a_1, \ldots, a_{h-u} \}$ and $\{
b_1, \ldots, b_{h-u} \}$.  Similarly, there is a bijection between
the sets $\{ a'_1, \ldots, a'_{h} \}$ and $\{ b'_1, \ldots, b'_{h}
\}$, and so the double representations~(\ref{ms:uniqueu})
and~(\ref{ms:uniquew}) are equivalent.  Thus, for every positive
integer $u$ there is at most one proper double representation of
the form~(\ref{ms:uniqueu}).

If $u<w,$ we obtain the double representation
\[
a_1+\cdots + a_{h-u}+ b'_1 + b'_2 + \cdots + b'_h = (w-u)x + a'_1 +
a'_2 + \cdots + a'_h + b_1+\cdots + b_{h-w}.
\]
Cancelling elements that appear on both sides of this equation, we
obtain a proper double representation of the form
\[
(w-u)x + a_1+\cdots + a_{h-w+u} = a'_1 + a'_2 + \cdots + a'_h,
\]
where $w-u \geq 1$ and $\{ a_1, \ldots, a_{h-w+u}, a'_1, a'_2,
\ldots, a'_h \} \subseteq A.$  We call this process the
``subtraction algorithm.''

Let $u$ be the smallest positive integer for which there exists a
proper double representation of the form~(\ref{ms:uniqueu}). Suppose that
there is a proper double representation of the
form~(\ref{ms:uniquew}) for some integer $w > u.$
By the division algorithm, we write $w = qu+r,$
where $0 \leq r < u.$   If $r \geq 1,$ then iteration of the subtraction algorithm above
yields a proper double representation in which the element $x$
appears exactly $r$ times, which contradicts the minimality of $u$.
It follows that $u$ must divide $w$.
Moreover, if there exists a proper double representation for some $w > u,$
then the subtraction algorithm produces a double representation with $w = 2u.$
Thus we have proper double representations of the form
\beq  \label{ms:addw}
ux + a_1+\cdots + a_{h-u} = a'_1 + a'_2 + \cdots + a'_h
\eeq
and
\beq  \label{ms:add2w}
2ux + b_1+\cdots + b_{h-2u} = b'_1 + b'_2 + \cdots + b'_h,
\eeq
where
\[
\{a_1,\ldots, a_{h-u}\} \cap \{ a'_1 , a'_2 , \ldots , a'_h \} = \emptyset
\]
and
\[
\{b_1,\ldots, b_{h-2u}\} \cap \{ b'_1 , b'_2 , \ldots , b'_h \} = \emptyset.
\]
Adding equations~(\ref{ms:addw}) and~(\ref{ms:add2w}) and
cancelling $ux,$ we obtain the following double representation of
length $2h-u$:
\[
ux + a'_1 + a'_2 + \cdots + a'_h + b_1+\cdots + b_{h-2u} =
a_1+\cdots + a_{h-u} + b'_1 + b'_2 + \cdots + b'_h.
\]
After subtracting $h-u$ equal terms on both sides of this equation, we must obtain
the proper double representation~(\ref{ms:addw}).  This means that on
the left side we must have $a_1 + a_2 + \cdots + a_{h-u}$.  Since
\[
\{a_1,\ldots, a_{h-u}\} \cap \{ a'_1 , a'_2 , \ldots , a'_h \} = \emptyset,
\]
it follows that the sequence $(a_1 , a_2 , \ldots , a_{h-u})$
is equivalent to the sequence $(b_1 , b_2 , \ldots , b_{h-2u})$,
which is absurd since $h-2u < h-u.$
This completes the proof.
\end{proof}

\bl \label{ms:lemma:2h-1}     Let $h \geq 2,$ and let $A$ be a
maximal $B_h$-subset of the $B_{2h-1,h-1}$-set $X$.  Let $x \in
X\setminus A,$ and let \beq \label{ms:uuu} ux + a_1+\cdots + a_{h-u}
= a'_1 + a'_2 + \cdots + a'_h \eeq be the unique proper double
representation of length $h$ with elements in $A \cup \{x\}.$  For
every
\[
a^{\ast} \in \{
a_1,a_2,\ldots,a_{h-u}, a'_1,a'_2,\ldots,a'_h \},
\]
the set
\[
\left( A \cup \{ x \} \right) \setminus \{a^{\ast} \}
\]
is a $B_h$-set contained in $X$. \el

\begin{proof}
If $\left( A \cup \{ x \} \right) \setminus \{a^{\ast} \} $ is not a
$B_h$-set, then there must exist a positive integer $v$ and elements
$b_1, \ldots, b_v,b'_1,\ldots,b'_h \in A \setminus \{a^{\ast}\} $
such that \beq  \label{ms:vvv} vx + b_1+\cdots + b_{h-v} = b'_1 +
b'_2 + \cdots + b'_h \eeq is a proper double representation in $X$.
Then (\ref{ms:uuu}) and~(\ref{ms:vvv}) are different proper double
representations of length $h$ in $X$, which contradicts
Lemma~\ref{ms:lemma:proper-one}.
\end{proof}

\bt \label{ms:theorem:bases} Let $h \geq 2,$ and let $X$ be a finite
$B_{2h-1,h-1}$-set contained in the abelian group $\Gamma.$ Then the
maximal $B_h$-subsets of $X$ have the same cardinality. \et

\begin{proof}
Let $\mathcal{M}_h(X)$ be the set of maximal $B_h$-sets contained in
$X$, and let
\[
m = \max\{ |C| : C \in \mathcal{M}_h(X) \}.
\]
We must prove that $|C| = m$ for every $C \in \mathcal{M}_h(X).$

Let $C \in \mathcal{M}_h(X),$ and let $C^{\ast}$ be the largest
subset of $C$ that is contained in a $B_h$-set $A$ of cardinality
$m$.  If $C^{\ast} = C,$ then the maximality of $C$ implies that $C
= A$, and so $|C| = m.$  If $C^{\ast} \neq C,$ then there exists $s
\in C \setminus A.$  By the maximality of $A$, the set $A \cup
\{s\}$ is not a $B_h$-set, and there exists a proper double
representation of the form
\[
ws+a_1+\cdots + a_{h-w} = a'_1 + \cdots + a'_h,
\]
where
\[
\{s, a_1,a_2,\ldots,a_{h-u} \} \cap  \{ a'_1,a'_2,\ldots,a'_h \} =
\emptyset.
\]
If \beq  \label{ms:badrep} \{a_1,a_2,\ldots,a_{h-u},
a'_1,a'_2,\ldots,a'_h \} \subseteq C^{\ast} \subseteq A, \eeq
then~(\ref{ms:badrep}) is a proper double representation of length
$h$ with elements in $C$, which contradicts the fact that $C$ is a
$B_h$-set.  Therefore, there exists an element
\[
a^{\ast} \in \{a_1,a_2,\ldots,a_{h-u}, a'_1,a'_2,\ldots,a'_h \}
\setminus C^{\ast}.
\]
By Lemma~\ref{ms:lemma:2h-1}, the set
\[
\left( A \cup \{ s \} \right) \setminus \{a^{\ast} \}
\]
is a $B_h$-set contained in $X$, and
\[
C^{\ast} \cup \{ s \} \subseteq \left( A \cup \{ s \} \right)
\setminus \{a^{\ast} \}.
\]
This is impossible, since
\[
C^{\ast} \cup \{ s \} \subseteq C,
\]
\[
|C^{\ast} \cup \{ s \}| = |C^{\ast}| + 1,
\]
and
\[
\left| \left( A \cup \{ s \} \right) \setminus \{a^{\ast} \} \right|
= |A| = m.
\]
Therefore, $C^{\ast} = C.$  This completes the proof.
\end{proof}

\section{Matroids of $B_h$-sets}
A {\em matroid} $M = M(X,\mathcal{I})$ consists of a finite set $X$
and a collection $\mathcal{I}$ of subsets of $X$ that satisfy the
following properties: \benum
\item[(i)] $\emptyset \in \mathcal{I},$
\item[(ii)]
If $B \in \mathcal{I}$ and $A \subseteq B,$ then $A \in
\mathcal{I},$
\item[(iii)]
If $A,B \in \mathcal{I}$ and $|A| < |B|,$ then there exists $b \in
B\setminus A$ such that $A \cup \{b\} \in \mathcal{I}.$ \eenum The
members of $\mathcal{I}$ are called the {\em independent sets} in
$X$. A {\em basis} for $X$ is a maximal independent set.
Condition~(iii) implies that all bases have the same cardinality.
The {\em rank} of the matroid $M$ is the cardinality of a basis for
$M$.

\bt \label{ms:theorem:matroid} Let $h \geq 2,$ and let $X$ be a
finite $B_{2h-1,h-1}$-subset of an abelian group. Let $\mathcal{I}$
be the collection of $B_h$-sets contained in $X$.   Then $M =
M(X,\mathcal{I})$ is a matroid. \et

\begin{proof}
Every subset of a $B_h$-set is a $B_h$-set, and the empty set is
also a $B_h$-set.  We must show that if $A$ and $B$ are
$B_h$-subsets of $X$ with $|A| < |B|,$  then there exists $b \in
B\setminus A$ such that $A \cup \{ b \}$ is a $B_h$-set.

Let
\[
X' = A \cup B.
\]
Then $X'$ is a $B_{2h-1,h-1}$-subset of $X$.  Let $m$ be the
cardinality of the maximal $B_h$-subsets of $X'$. Let $A^{\ast}$ be
a maximal $B_h$-subset of $X'$ that contains $A$. Then
\[
|A| < |B| \leq m = |A^{\ast}|,
\]
and so there exists an element
\[
b \in A^{\ast} \setminus A  \subseteq X'\setminus A = B \setminus A.
\]
Then $A \cup \{ b\} \subseteq A^{\ast},$ and so $A \cup \{ b\}$ is a
$B_h$-set. This completes the proof.
\end{proof}

Let $M = M(X,\mathcal{I})$ be a matroid.  For every positive integer
$k$, let $\mathcal{I}^{(k)}$ be the set of all unions of $k$
independent subsets of $X$, that is, all sets of the form $I_1 \cup
I_2 \cup \cdots \cup I_k,$ where $I_1, I_2,\ldots, I_k \in
\mathcal{I}.$ Then $M^{(k)} = M(X,\mathcal{I}^{(k)})$ is also a
matroid on the set $X$ (Welsh~\cite[Section 8.3]{wels76}). We denote
the rank of the matroid $M^{(k)}$ by $\rho_k.$   Then $\rho_k$ is
the cardinality of the largest subset of $X$ that can be written as
the union of $k$ independent sets in $X$.

The {\em covering number} of a set $S$ contained in $X$ is the
smallest integer $k$ such that $S$ can be written as the union of
$k$ independent subsets of $X$.  If $\{x\} \in \mathcal{I}$ for
every $x \in X,$ then the covering number exists, and the covering
number of $S$ is at most $|S|.$  The set $S$ has covering number $k$
if and only if $k$ is the smallest integer such that $S$ is an
independent set in the matroid $M^{(k)} $.  The set $X$ has covering
number $k$ if and only if
\[
\rho_1 < \rho_2 < \cdots < \rho_k = |X|.
\]

Let $ X $ be a $B_{2h-1,h-1} $-set contained in an abelian group.
For every subset $S$ of $X$, we define the {\em $B_h$-covering
number} of $S$ as the smallest integer $k$ such that $S = A_1 \cup
\cdots \cup A_k,$ where $A_1,\ldots, A_k$ are $B_h$-sets.  Since
$\{x\}$ is a $B_h$-set for all $x \in X,$ it follows that every
subset of $X$ has a finite $B_h$-covering number.

\bt Let $X$ be a $B_{2h-1,h-1}$-set contained in an abelian group,
and let $\ell$ be the $B_h$-covering number of $X$.  For every
positive integer $k \leq \ell$ there is a number $n_X(k)$ such
that if $S$ is a maximal subset of $X$ whose $B_h$-covering number
is $k$, then $|S|= n_X(k).$ \et

\begin{proof}
By Theorem~\ref{ms:theorem:matroid}, $M = M(X,\mathcal{I})$ is a
matroid, where $\mathcal{I}$ is the set of $B_h$-subsets of $X$. Let
$n_X(k)$ denote the cardinality of the bases in the matroid
$M^{(k)}.$  The maximal subsets of $X$ with $B_h$-covering number
$k$ are precisely the bases in the matroid $M^{(k)}.$  This
completes the proof.
\end{proof}

Let $I_1, I_2,\ldots, I_k$ be independent sets in a matroid $M =
M(X,\mathcal{I})$.  We define $I'_1 = I_1$ and $I'_j = I_j
\setminus (I_1 \cup \cdots \cup I_{j-1})$ for $j = 2,\ldots,k.$
Since every subset of an independent set is independent, it
follows that the sets $I'_1, I'_2,\ldots, I'_k$ are pairwise
disjoint independent sets in $M$, and  $I_1 \cup I_2 \cup \cdots
\cup I_k = I'_1 \cup I'_2 \cup \cdots \cup I'_k$. Therefore, every
independent set in the matroid $M^{(k)}$ can be written as the
union of $k$ pairwise disjoint independent sets in $M$.  In
particular, if $X$ has covering number $k,$ then $X$ is the union
of $k$ pairwise disjoint nonempty independent subsets of $X$.

Let $\mu = (\mu_1,\ldots, \mu_r)$ be a partition of $|X|,$ that is,
$\mu_1,\mu_2,\ldots, \mu_r$ are positive integers such that $\mu_1 +
\mu_2 +  \cdots + \mu_r = |X|$ and $\mu_1 \geq \mu_2 \geq \cdots
\geq \mu_r.$   A {\em $\mu$-covering} of the matroid $M =
M(X,\mathcal{I})$ consists of $r$ pairwise disjoint independent sets
$I_1, I_2,\ldots, I_r$ such that
\[
X = I_1 \cup I_2 \cup \cdots \cup I_r
\]
and
\[
|I_j| = \mu_j \qquad\text{for $j = 1,2,\ldots,r.$}
\]
Let $k$ be the covering number of the matroid $M$.   Dias da
Silva~\cite{dias90} proved that there exists a $\mu$-covering of $X$
if and only if $k \leq r$ and $\rho_j \geq \mu_1 + \cdots + \mu_j$
for $j = 1,\ldots, k.$

\bt Let $X$ be a $B_{2h-1,h-1}$-set contained in an abelian group,
and let $k$ be the $B_h$-covering number of $X$.  For $j =
1,\ldots,k,$ let $\rho_j$ denote the maximum cardinality of a union
of $j$ $B_h$-subsets of $X$.  Let $\mu = (\mu_1,\ldots, \mu_r)$ be
any partition of $|X|.$ There exist pairwise disjoint $B_h$-sets
$I_1,\ldots,I_r$ such that
\[
X= I_1 \cup \cdots \cup I_r
\]
and $|I_j| = \mu_j$ for $j = 1,\ldots, r$ if and only if $r \geq
k$ and $\rho_j \geq \mu_1 + \cdots + \mu_j$ for $j = 1,\ldots, k.$
\et

\begin{proof}
This follows immediately from the fact that the $B_h$-sets are the
independent sets of a matroid on $X$.
\end{proof}

\providecommand{\bysame}{\leavevmode\hbox
to3em{\hrulefill}\thinspace}
\providecommand{\MR}{\relax\ifhmode\unskip\space\fi MR }
\providecommand{\MRhref}[2]{%
  \href{http://www.ams.org/mathscinet-getitem?mr=#1}{#2}
} \providecommand{\href}[2]{#2}

\end{document}